\documentclass[12pt,oneside,english]{amsart}
\usepackage[latin1]{inputenc}
\usepackage{amssymb}

\makeatletter
\newtheorem{theorem}{Theorem}
\newtheorem{lemma}{Lemma}
\newtheorem{example}{Example}
\newtheorem{definition}{Definition}
\newtheorem{condition}{Condition}

\newtheorem{postulate}{Postulate}
\newtheorem{criterion}{Criterion}
\newtheorem{remark}{Remark}
\newtheorem{conjecture}{Conjecture}

\newtheorem{corollary}{Corollary}

\numberwithin{equation}{section}
\usepackage{babel}

\makeatother
\begin{document}
\baselineskip=17pt

\title[Theorems on twin primes\enskip-\enskip dual case]{Theorems on twin primes\enskip-\enskip dual case }

\author{Vladimir Shevelev}
\address{Department of Mathematics \\Ben-Gurion University of the
 Negev\\Beer-Sheva 84105, Israel. e-mail:shevelev@bgu.ac.il}

\subjclass{Primary 11A41, secondary 11B05}

\begin{abstract}
We prove dual theorems to theorems proved by author in \cite {5}. Beginning with Section 10,
we introduce and study so-called "twin numbers of the second kind" and a postulate for them.
We give two proofs of the infinity of these numbers and a sufficient condition for
truth of the postulate; also we pose several other conjectures. Finally, we consider
a conception of axiom of type "AiB".

\end{abstract}
\maketitle

\section{Introduction }

In [3] we posed, in particular, the following conjecture
\begin{conjecture}\label{1}
Let $\tilde{c}(1)=2$ and for $n\geq 2,$
$$\tilde{c}(n)=\tilde{c}(n-1)+\begin{cases}\gcd (n, \enskip \tilde{c}(n-1)),\enskip if \;\; n\enskip is\enskip odd
\\\gcd (n-2,\enskip \tilde{c}(n-1)),\enskip if\;\; n\enskip is\enskip even.\end{cases} $$
 Then every record (more than 3) of the values of difference $\tilde{c}(n)-\tilde{c}(n-1)$ is greater of twin primes.
 \end{conjecture}
 The first records are (cf. sequence A167495 in [6])
 \begin{equation}\label{1.1}
  5, 13, 31, 61, 139, 283, 571, 1153, 2311, 4651, 9343, 19141, 38569,...
  \end{equation}
  We use the same way as in our paper [5] which is devoted to study a sequence dual to the now considered one.
  Our observations of the behavior of sequence $\{\tilde{c}(n)\}$ are the following:\newline\newline

  1) In some sequence of arguments $\{m_i\}$ we have $\frac {\tilde{c}(m_i)-3} {m_i-3}=3/2.$ These values of arguments we call \slshape the fundamental points. \upshape The first fundamental point are
  $$ 7, 27, 63, 123, 279, 567, 1143, 2307, 4623, 9303, 18687,...$$
  2)For every two adjacent fundamental points $m_j<m_{j+1},$ we have $m_{j+1}\geq 2m_j-3.$\newline
  3) For $i\geq2,$ the numbers $\frac {m_i-5} {2},\enskip\frac {m_i-1} {2}$ are twin primes (and, consequently, $m_i\equiv3\pmod {12}$).\newline
  4) In points $m_i+1$ we have $\tilde{c}(m_i+1)-\tilde{c}(m_i)=\frac{m_i-1} {2}.$ These increments we call \slshape the main increments\upshape \enskip of sequence $\{\tilde{c}(n)\},$ while other nontrivial (i.e.more than 1) increments we call \slshape the minor increments.\upshape\newline
  5)For $i\geq2,$ denote $h_i$ the number of minor increments between adjacent fundamental points $m_i$ and $m_{i+1}$ and $T_i$ the sum of these increments. Then $T_i\equiv h_i\pmod6.$\newline
  6) For $i\geq2,$ the minor increments between adjacent fundamental points $m_i$ and $m_{i+1}$ could occur only before $m_{i+1}-\sqrt{2( m_{i+1}-1)}-2.$ \upshape\newline\newline
  \indent Below we show that the validity of all these observations follow only from 6).
  \begin{theorem}\label{1}
If observation 6) is true then observation 1)-5) are true as well.
 \end{theorem}
 \begin{corollary}\label{1}
 If $1)$ observation 6) is true and $2)$ the sequence $\{\tilde{c}(n)\}$ contains infinitely many fundamental points, then there exist infinitely many twin primes.
\end{corollary}
 Besides, in connection with Conjecture 1 we think that
 \begin{conjecture}\label{2}
  For $n\geq16,$ the main and only main increments are the record differences $\tilde{c}(n)-\tilde{c}(n-1).$
 \end{conjecture}

\section{Proof of Theorem 1}
 We use induction. Suppose $n_1\geq28$ is a number of the form 12l+4 (for $n_1<28$ the all observations are verified directly).
Let $n_1-1$ is a fundamental point and for $n:=\frac{n_1-4} {2},\enskip n\mp1$ are twin primes. Thus $$\tilde{c}(n_1-1)=\frac{3} {2}(n_1-4)+3=\frac {3} {2}n_1-3.$$
Since $n_1$ is even and
$$\gcd(\frac {3} {2}n_1-3, n_1-2)=\frac {n_1} {2}-1,$$
 then we have a main increment such that
 \begin{equation}\label{2.1}
\tilde{c}(n_1)=2n_1-4.
 \end{equation}

Here we distinguish two cases:\newpage
\bfseries A \mdseries) Up to the following fundamental point there are only trivial increments. The inductive step in this case we formulate as the following.
 \begin{theorem}\label{2}
If $27\leq m_i<m_{i+1}$ are adjacent fundamental points without miner increments between them, then
i)\enskip$m_{i+1}=2m_i-3;$\newline
ii)\enskip If  \enskip$\frac{m_i-5} {2},\enskip\frac{m_i-1} {2}$ are twin primes, then $\frac {m_{i+1}-5} {2},\enskip\frac{m_{i+1}-1} {2}$ are twin primes as well.
 \end{theorem}
 Note that really, for the first time, Case \bfseries A \mdseries) appears for $m_3=63,$ such that, by Theorem 2, we have two pairs of twin primes: (29,31), (59,61).

\bfseries Inductive step in case A \mdseries)\newline Continuing (2.1), we have
$$\tilde{c}(n_1+1)=2n_1-3,$$
$$\tilde{c}(n_1+2)=2n_1-2,$$
$$...$$
$$\tilde{c}(2n_1-5)=3n_1-9,$$
Since $\frac {3n_1-12} {2n_1-8}=3/2,$ then we conclude that $2n-1-5$ is the second fundamental point in the inductive step. By the definition of the sequence, denoting $n_2=2n_1-4,$ we have
\begin{equation}\label{2.2}
\tilde{c}(n_2)=2n_2-4.
\end{equation}
Note that, since $n_1=12l+4,$ then $n_2=12l_1+4,$ where $l_1=2l.$\newline
\indent Furthermore, from the run of formulas (2.2) we find for $3\leq j\leq\frac{n_1-2} {2}:$
$$\tilde{c}(2n_1-2j-1)=3n_1-2j-5,$$
$$\tilde{c}(2n_1-2j)=3n_1-2j-4.$$
This means that
$$\gcd(2n_1-2j-2,\enskip 3n_1-2j-5)=1,\enskip i.e.\enskip \gcd(j-2,\enskip n_1-3)=1.  $$
Note that, for the considered values of $n_1$ we have $\frac{n_1-2} {2}\geq\sqrt{n_1-3},$ then $n_1-3=\frac {n_2-2} {2}$ is prime.\newline
\indent On the other hand,
$$\tilde{c}(2n_1-2j)=3n_1-2j-4,$$
$$\tilde{c}(2n_1-2j+1)=3n_1-2j-3.$$
Thus, for $7\leq j\leq\frac{n_1-2} {2},$
$$\gcd(2n_1-2j+1,\enskip 3n_1-2j-4)=1,\enskip i.e.\enskip \gcd(2j-11,\enskip n_1-5)=1.$$

Here, for the considered values of $n_1$ we also have $2n_1-13\geq\sqrt{n_1-5},$ then $n_1-5=\frac {n_2-6} {2}$ is prime. $\blacksquare$
\newpage

 \bfseries B \mdseries) Up to the following fundamental point we have some minor increments.

   The inductive step we formulate as following.
   \begin{theorem}\label{3}
Let observation 6) be true. If $7\leq m_i<m_{i+1}$ are adjacent fundamental points with a finite number of minor increments between them, then \newline
i)\enskip$m_{i+1}\geq2m_i;$\newline
ii)\enskip If  \enskip$\frac{m_i-5} {2},\enskip\frac{m_i-1} {2}$ are twin primes, then $\frac {m_{i+1}-5} {2},\enskip\frac{m_{i+1}-1} {2}$ are twin primes as well.
 \end{theorem}
 Thus the observation 2) will be proved in frameworks of the induction.\newline

 \bfseries Inductive step in case B \mdseries)

 Let in the points $n_1+l_j\enskip j=1,...,h,$ before the second fundamental point we have the minor increments $t_j,\enskip j=1,...,h.$ We have ( starting with the first fundamental point $n_1-1$)

$$\tilde{c}(n_1-1)=\frac {3} {2}n_1-3,$$
$$\tilde{c}(n_1)=2n_1-4,$$
$$\tilde{c}(n_1+1)=2n_1-3,$$
$$...$$

$$\tilde{c}(n_1+l_1-1)=2n_1+l_1-5.$$

\begin{equation}\label{2.3}
\tilde{c}(n_1+l_1)=2n_1+l_1+t_1-5,
\end{equation}

$$\tilde{c}(n_1+l_1+1)=2n_1+l_1+t_1-4,$$

$$...$$

$$\tilde{c}(n_1+l_2-1)=2n_1+l_2+t_1-6,$$

\begin{equation}\label{2.4}
\tilde{c}(n_1+l_2)=2n_1+l_2+t_1+t_2-6,
\end{equation}
$$...$$
$$\tilde{c}(n_1+l_h-1)=2n_1+l_h+t_1+...+t_{h-1}-h-5,$$\newpage
\begin{equation}\label{2.5}
\tilde{c}(n_1+l_h)=2n_1+l_h+t_1+...+t_h-h-4,
\end{equation}
$$\tilde{c}(n_1+l_h+1)=2n_1+l_h+t_1+...+t_h-h-3,$$
$$...$$
\begin{equation}\label{2.6}
\tilde{c}(2n_1+2T_h-2h-5)=3n_1+3T_h-3h-9,
\end{equation}
where
\begin{equation}\label{2.7}
T_h=t_1+...+t_h.
\end{equation}
It is easy to see that $2n_1+2T_h-2h-5$ is the second fundamental point in the inductive step.
Furthermore, subtracting 2 from the even number $2n_1+2T_h-2h-4,$ we see that
$$\gcd(2n_1+2T_h-2h-6,\enskip 3n_1+3T_h-3h-9)=n_1+T_h-h-3.$$
 Thus in the point $n_2:=2n_1+2T_h-2h-4$ we have the second main increment (in framework of the inductive step):
\begin{equation}\label{2.8}
\tilde{c}(2n_1+2T_h-2h-4)=4n_1+4T_h-4h-12.
\end{equation}
   Note that, for $n\geq2,$ we have $\tilde{c}(n)\equiv n\pmod2.$ Therefore, $T_h\geq3h$ and for the second fundamental point $n_2-1=2n_1+2T_h-2h-5$ we find
 \begin{equation}\label{2.9}
n_2-1\geq2(n_1-1)+4h-3.
\end{equation}
This in frameworks of the induction confirms observation 2).\newline\newline
\indent Now, in order to finish the induction, we prove the primality of numbers $\frac {n_2-6} {2}=n_1+T_h-h-5$ and
$\frac {n_2-2} {2}=n_1+T_h-h-3.$\newline
\indent From the run of formulas (2.5)-(2.6) for $7\leq j\leq\frac{n_1+2T_h-2h-l_h} {2}$
(we cannot cross the upper boundary of the last miner increment) we find

$$\tilde{c}(2n_1+2T_h-2h-2j)=3n_1+3T_h-3h-2j-4,$$
$$\tilde{c}(2n_1+2T_h-2h-2j+1)=3n_1+3T_h-3h-2j-3.$$
Thus, for $7\leq j\leq\frac{n_1+2T_h-2h-l_h} {2},$
$$\gcd(2n_1+2T_h-2h-2j+1,\enskip 3n_1+3T_h-3h-2j-4)=1,$$ i.e.$$\gcd(2j-11,\enskip n_1+T_h-h-5)=1.$$
For the most possible $j=\frac{n_1+2T_h-2h-l_h-1} {2}$ (it is sufficient to consider the case of odd $l_h$) we should have
$$2j-11=n_1+2T_h-2h-l_h-12\geq\sqrt{n_1+T_h-h-5},$$
\newpage or, since $n_2=2n_1+2T_h-2h-4,$ then we should have $n_2-n_1-l_h-8\geq\sqrt{\frac{n_2-6} {2}},$
 i.e.
\begin{equation}\label{2.10}
n_1+l_h\leq n_2-\sqrt{\frac{n_2-6} {2}}-8,
\end{equation}
 Since $n_2\geq28,$ then this condition, evidently, follows from observation 6) which is written in terms of the fundamental points $m_i=n_i-1.$
Thus from observation 6) we indeed obtain the primality of $\frac{n_2-6} {2}=n_1+T_h-h-5.$ \newline
Furthermore,
$$\tilde{c}(2n_1+2T_h-2h-2j+1)=3n_1+3T_h-3h-2j-3,$$
$$\tilde{c}(2n_1+2T_h-2h-2j+2)=3n_1+3T_h-3h-2j-2.$$
Thus, for $6\leq j\leq\frac{n_1+2T_h-2h-l_h} {2},$
$$\gcd(2n_1+2T_h-2h-2j,\enskip 3n_1+3T_h-3h-2j-3)=1,$$ i.e.$$\gcd(j-3,\enskip n_1+T_h-h-3)=1.$$
For the most possible $j=\frac{n_1+2T_h-2h-l_h-1} {2}$ (here again sufficiently to consider the case of odd $l_h$) we should have
 $$\frac{n_1+2T_h-2h-l_h-1} {2}-3\geq\sqrt{\frac{n_2-2} {2}},$$
 or
\begin{equation}\label{2.11}
n_1+l_h\leq n_2-\sqrt{{2(n_2-2)}}-3.
\end{equation}
This coincides with observation 6). Thus $\frac{n_2-2} {2}$ is prime as well. This completes proof of Theorem 1 $\blacksquare$\newline
\indent Note that in [5] we used the Rowland method [2] to obtain an independent from observation 6) proof of the primality of the greater number. Here we give a parallel proofs for both of numbers.
\begin{corollary}\label{2}
 If $p_1<p_2$ are consecutive seconds of twin primes giving by Theorem 1, then $p_2\geq 2p_1-1.$
\end{corollary}
\bfseries Proof.  \mdseries The corollary easily follows from (2.9).$\blacksquare$
\begin{corollary}
 $$T_h\equiv h\pmod6.$$
 \end{corollary}\label{3}
\bfseries Proof.  \mdseries  The corollary follows from the well known fact that the half-sum of twin primes not less than 5 is a multiple of 6. Therefore,  $n_1+T_h-h-4\equiv0\pmod6.$ Since, by the condition, $n_1\equiv4\pmod{12},$ then we obtain the corollary.$\blacksquare$\newpage
\indent Now the observation 5) follows in the frameworks of the induction. The same we can say about observation 4).
\indent The observed weak excesses of the exact estimate of Corollary 2 indicate to the smallness of $T_h$ and confirm, by Theorem 1, Conjecture 1.

\section{A rule for constructing a pair of twin primes $p,\enskip p+2$ by a given integer $m\geq4$ such that $p+2\geq m$}
One can consider a simple rule for constructing a pair of twin primes $p,\enskip p+2$ by a given integer $m\geq4$ such that $p+2\geq m$ quite similar to one over sequence $\{c(n)\}$ (see Section 6 in [5]). To this aim, with $m$ we associate the sequence

$$\tilde{c}^{(m)}(1)=m;\enskip for\enskip n\geq2,$$
 \begin{equation}\label{3.1}
\tilde{c}^{(m)}(n)=\tilde{c}^{(m)}(n-1)+\begin{cases}\gcd (n, \enskip \tilde{c}^{(m)}(n-1)),\enskip if \;\; n\enskip is\enskip even
\\\gcd (n-2,\enskip \tilde{c}^{(m)}(n-1)),\enskip if\;\; n\enskip is\enskip odd.\end{cases}
\end{equation}
Thus for every $m$ this sequence has the the same formula that the considered one but with another initial condition.
Our observation is the following.
\begin{conjecture}\label{3}
Let $n^*,$ where $n^*=n^*(m),$ be point of the last nontrivial increment of $\{\tilde{c}^{(m)}(n)\}$ on the set $A_m=\{1,...,m-3\}$ and $n^*=1,$ if there is not any nontrivial increment on $A_m.$
Then numbers $\tilde{c}^{(m)}(n^*)-n^*\mp1$ are twin primes.
 \end{conjecture}
Evidently, $c^{(m)}(n^*)-n^*+1\geq m$ and the equality holds if and only if $n^*=1.$\newline
 \indent The following examples show that, for the same $m,$ the pair of twin primes which is obtained by the considered rule, generally speaking, differs from one which is obtained by the corresponding rule in [5].

  \begin{example}\label{1}
Let $m=577.$ Then $n^*=51$ and $\tilde{c}^{(m)}(n^*)=669.$ Thus numbers $669-51\mp1$ are twin primes $(617,\enskip619),$ while by the rule in $[5]$ we had another pair: $(881,\enskip883).$
 \end{example}

  \begin{example}\label{2}
Let $m=3111.$ Then $n^*=123$ and $\tilde{c}^{(m)}(n^*)=3513.$ Thus numbers $3513-123\mp1$ are twin primes $(3389,\enskip3391),$ while by the rule in $[5]$ we have another pair: $(3119,\enskip3121).$
 \end{example}
 The case of $n^*=1$ we formulate as the following criterion, which is proved quite similar to Criterion 1 [5].\newpage
 \begin{criterion}\label{1}A positive integer $m>3$ is a greater of twin primes if and only if all the points $1,...,m-3$ are points of trivial increments of sequence $\{\tilde{c}^{(m)}(n)\}.$
\end{criterion}
\section{A new sequence and an astonishing observation}
Consider the sequence which is defined by the recursion:
$$f(1)=2 \enskip and,\enskip for \enskip n\geq2,$$
$$f(n)=f(n-1)+\begin{cases}\gcd (n, \enskip f(n-1)+2),\enskip if \;\; n\enskip is\enskip even
\\\gcd (n,\enskip f(n-1)),\enskip if\;\; n\enskip is\enskip odd.\end{cases} $$

Here the even points $m_i\neq8$ in which $f(m_i)/m_i=3/2$ we call the fundamental points. The increments $\frac{m_i+2} {2}$ in the points $n_i=m_i+2$ are called main increments and other nontrivial (i.e. different from 1) increments we call miner increments. This sequence also could be studied by method of [5]. It is easy to verify that the nontrivial increments of this sequence differs from ones of the above considered sequence $\{\tilde{c}(n)\}.$ But, our  observations show that a very astonishing fact,probably, is true: \slshape all records more than 7 for sequences $\{\tilde{c}(n)\}$ and $\{f(n)\}$ coincide! \upshape We think that it is a deep open problem.\newline\newline

\section{Some other new sequences  connected with twin primes}
Here we present three additional new sequences of the considered type, the records of which are undoubtedly connected with twin primes.\newline
1)$$g(1)=2 \enskip and,\enskip for \enskip n\geq2,$$
$$g(n)=g(n-1)+\begin{cases}\gcd (n, \enskip g(n-1)+2),\enskip if \;\; n\enskip is\enskip even
\\\gcd (n-2,\enskip g(n-1)+2),\enskip if\;\; n\enskip is\enskip odd.\end{cases} $$

2)$$h(1)=2 \enskip and,\enskip for \enskip n\geq2,$$
$$h(n)=h(n-1)+\begin{cases}\gcd (n-2, \enskip h(n-1)+2),\enskip if \;\; n\enskip is\enskip even
\\\gcd (n,\enskip h(n-1)+2),\enskip if\;\; n\enskip is\enskip odd.\end{cases} $$
3)$$i(1)=2 \enskip and,\enskip for \enskip n\geq2,$$
$$i(n)=i(n-1)+\gcd (n, \enskip i(n-1)+2(-1)^n). $$
Note that, all records of the second sequence are, probably, the firsts of twin primes.
\newpage
\section{A theorem on twin primes which is independent on observation of type 6)}
Here we present a new sequence $\{\tilde{a}(n)\}$ with the quite analogous definition of fundamental and miner points for which Corollary 1 is true in a stronger formulation. Using a construction close to those ones that we considered in [4], consider the sequence defined as the following:
 $\tilde{a}(39)=57$ and for $n\geq 23,$
\begin{equation}\label{6.1}
\tilde{a}(n)=\begin{cases}\tilde{a}(n-1)+1,\enskip if \;\;\gcd(n-(-1)^n-1,\enskip \tilde{a}(n-1))=1;
\\2(n-2)\;\; otherwise\end{cases}.
\end{equation}
The sequence has the following first nontrivial differences
$$19, 6, 2, 43,5,2,2,7,6,2,103,5,2,2,18,2,229,6,2,463,...$$

\begin{definition} A point $m_i$ is called \upshape a fundamental point of sequence (6.1),\slshape \enskip if it has the form $m_i=12t+3$  \enskip$and$ \enskip $\tilde{a}(m_i)-3=\frac{3}{2}(m_i-3).$ The increments in the points $m_i+1$ we call the \upshape main increments.\slshape \enskip Other nontrivial increments we call \upshape miner increments.
\end{definition}
The first two fundamental points of sequence (6.1) are 39 and 87.

\begin{theorem}\label{4}
If the sequence $\{{\tilde{a}}(n)\}$ contains infinitely many fundamental points, then there exist infinitely many twin primes.
 \end{theorem}
 \bfseries Proof.  \mdseries We use induction. Suppose, for some $i\geq1,$ the numbers $\frac {m_i-3} {2}\mp1$ are twin primes.
 Put $n_i=m_i+1.$ Then $n_i\equiv4\pmod{12}$ and we have

$$\tilde{a}(n_i-1)=\frac{3} {2}n_i-3,$$
$$\tilde{a}(n_i)=2n_i-4,$$
We see that the main increment is $\frac{n_i-2}{2}.$ By the condition, before $m_{i+1}$ we can have only a finite set if miner increments. Suppose that, they are in the points $n_i+l_j, j=1,...,h_i.$ Then, by (6.1), we have
$$\tilde{a}(n_i+1)=2n_i-3,$$
$$...$$
$$\tilde{a}(n_i+l_1-1)=2n_i+l_1-5,$$
$$\tilde{a}(n_i+l_1)=2n_i+2l_1-4,$$
$$...$$
$$\tilde{a}(n_i+l_2-1)=2n_i+l_1+l_2-5,$$
$$\tilde{a}(n_i+l_2)=2n_i+2l_2-4,$$
\newpage
$$...$$

$$\tilde{a}(n_i+l_h-1)=2n_i+l_{h-1}+l_h-5,$$
\begin{equation}\label{6.2}
\tilde{a}(n_i+l_h)=2n_i+2l_h-4,
\end{equation}

$$...$$

\begin{equation}\label{6.3}
\tilde{a}(n_{i+1}-1)=\frac{3}{2}n_{i+1}-3,
\end{equation}
\begin{equation}\label{6.4}
\tilde{a}(n_{i+1})=2n_{i+1}-4.
\end{equation}
Note that, in every step from (6.2) up to (6.3) we  add 1 simultaneously to values of the arguments and of the right hand sides.  Thus in the fundamental point $m_{i+1}=n_{i+1}-1$ we have
$$n_i+l_h+x=n_{i+1}-1$$
and
$$2n_i+2l_h-4+x=\frac{3}{2}n_{i+1}-3$$
such that
\begin{equation}\label{6.5}
n_{i+1}=2n_i+2l_h-4.
\end{equation}
Now we should prove that the numbers
$$\frac{m_{i+1}-3}{2}\mp1= \frac {n_{i+1}-4} {2}\mp1$$
i.e.
$$n_i+l_h-5, \enskip n_i+l_h-3$$
are twin primes. \newline
We have

$$\tilde{a}(n_i+l_h+t)=2n_i+2l_h-4+t,$$
\begin{equation}\label{6.6}
\tilde{a}(n_i+l_h+t+1)=2n_i+2l_h-3+t,
\end{equation}
where $0\leq t\leq n_i+l_h-7.$
Distinguish two case.\newline
1) Let $l_h$ be even. Then, for even values of $t$ the numbers $n_i+l_h+t+1$ are odd and from equalities (6.6) we have
$$\gcd(n_i+l_h+t+1,\enskip2n_i+2l_h-4+t)=1.$$
or
$$\gcd(n_i+l_h+t+1,\enskip n_i+l_h-2+t/2)=1$$\newpage
and
$$\gcd(t/2+3,\enskip n_i+l_h-5)=1,\enskip 0\leq t/2\leq (n_i+l_h-7)/2.$$
Thus $n_i+l_h-5$ is prime.\newline
On the other hand, for odd values of $t,$ taking into account that $n_i+l_h+t+1$ is even, from equalities (6.6) we have
$$\gcd(n_i+l_h+t-1,\enskip2n_i+2l_h-4+t)=1,$$

$$\gcd(2n_i+2l_h+2t-2,\enskip2n_i+2l_h-4+t)=1$$
and
$$\gcd(t+2,\enskip n_i+l_h-3)=1,\enskip 0\leq t\leq n_i+l_h-7,\enskip t\equiv1\pmod2.$$
Thus $n_i+l_h-3$ is prime as well and the numbers $n_i+l_h-5,\enskip n_i+l_h-3$ are indeed twin primes. \newline

2) Let $l_h$ be odd. Then, using again equalities (6.6), by the same way, we show that the numbers $n_i+l_h-5,\enskip n_i+l_h-3 $ are twin primes. \newline
\indent Besides, note that $n_i+l_h-4\equiv0\pmod6$ and, thus $m_{i+1}=n_{i+1}-1=2n_i+2l_h-5\equiv3 \pmod{12}.$ This completes the induction.$\blacksquare$
\section{Algorithm without trivial increments}
Sequences of the considered type in this paper and in [5] contain too many points of trivial $1$-increments. For example, 10000 terms of sequence $\{{\tilde{a}}(n)\}$ give only 8 pairs of twin primes. Therefore, the following problem is actual from the computation point of view just as from the research point of view : to accelerate this algorithm for receiving of twin primes by the omitting of the trivial increments. Below we solve this problem.
\begin{lemma}\label{1}
If sequence $\{{\tilde{a}}(n)\}$ has a miner increment $\Delta$ in even point, then $\Delta$ is prime.
\end{lemma}
\bfseries Proof.  \mdseries Let even N be a point of a miner increment and $M=N-k$ be a point of the previous nontrivial increment. We distinguish two cases: $M$ is even and  $M$ is odd.\newline
a)Let $M$ be even. Then we have
 $${\tilde{a}}(M)=2M-4,$$
 $${\tilde{a}}(M+1)=2M-3,$$
$$...$$
$${\tilde{a}}(M+k-1)=2M+k-5,$$\newpage
\begin{equation}\label{7.1}
{\tilde{a}}(N)={\tilde{a}}(M+k)=2M+2k-4,
\end{equation}
where $k$ is the least positive integer for which the point $M+k$ is the point of a nontrivial increment. We see that
$$\Delta=\Delta(N)=k+1.$$
Since in this case $k$ is even, then
$$\gcd(M+k-2,\enskip 2M+k-5)=d>1$$
and, therefore,
$$\gcd(k+1,\enskip M-3)=d>1.$$
Thus some prime divisor $P$ of $M-3$ divides $k+1$ and, therefore, $k+1\geq P.$ All the more,
$$k+1\geq p,$$
where $p$ is the least prime divisor of $M-3.$ Since in the considered case $M-3$ is odd, then $p$ is odd. But, since $p-2\leq k-1,$ then in the run of formulas (7.1) there is
the following
$${\tilde{a}}(M+p-2)=2M+p-6.$$
Nevertheless, the following value of argument is $M+p-1\equiv0\pmod2$ and both of the numbers $M+p-3$ and $2M+p-6$ are multiple of $p.$ This means that $k\leq p-1,$ such that we have
$$\Delta=\Delta(N)=k+1=p. $$
2) $M$ is odd. This case is considered quite analogously. Note that here $p\geq2.$
$\blacksquare$

\begin{lemma}\label{2}
 Let sequence $\{{\tilde{a}}(n)\}$ have a miner increment $\Delta$ in odd point. If the sequence has the previous nontrivial increment in even point, then $\Delta$ is even such that $(\Delta+4)/2$ is prime.
\end{lemma}
\bfseries Proof.  \mdseries Let odd N be a point of a miner increment and $M=N-k\equiv0\pmod2$ be a point of the previous nontrivial increment.
 Then we again have the run of formulas (7.1).
 Since here $k$ is odd, then
$$\gcd(M+k,\enskip 2M+k-5)=d>1$$
and, therefore,
$$\gcd((k+5)/2,\enskip M-5)=d>1$$

Thus some prime divisor $P$ of $M-5$ divides $(k+5)/2$ and, therefore, $k+5\geq 2P.$ All the more,\newpage
$$k+5\geq 2p,$$
where $p$ is the least prime divisor of $M-5.$ Since in the considered case $M-5$ is odd, then $p$ is odd. But in the run of formulas (7.1) there is the following
$${\tilde{a}}(M+2p-6)=2M+2p-10.$$
Nevertheless, the following value of argument is $M+2p-5\equiv0\pmod1$ and both of the numbers $M+2p-5$ and $2M+2p-10$ are multiple of $p.$ This means that $k\leq 2p-5,$ such that we have
$$\Delta(N)=k+1=2p-4. $$

$\blacksquare$\newline
\indent Quite analogously we obtain the following lemma.
\begin{lemma}\label{3}
 Let sequence $\{{\tilde{a}}(n)\}$ have a miner increment $\Delta$ in odd point. If the sequence has the previous nontrivial increment in odd point, then $\Delta$ is odd such that $\Delta+4$ is prime.
\end{lemma}
\begin{remark}
A little below we shall see that actually for nontrivial increments the conditions of Lemma 3 do not appear. But the proof
of Lemma 3 plays its role!
\end{remark}
Note now that in proofs of Lemmas 1-3 $p$ is always the least prime divisor of $M-5$ or $M-3,$ where $M$ is point of the "previous nontrivial increment," we obtain the following algorithm for the receiving of twin primes.
\begin{theorem}\label{5}
1) Let $n_m $ be point of the $m$-th main increment of sequence $\{{\tilde{a}}(n)\}$ and $P_m$ be the least prime divisor of the product $(n_m-5)(n_m-3).$ Then the first point $N_1$ of miner increment is
\begin{equation}\label{7.2}
N_{1}=\begin{cases}n_m+P_m-1,\enskip if \;\; \enskip P_m|(n_m-3),
\\n_m+2P_m-5,\enskip  if\;\; P_m|(n_m-5).\end{cases}
\end{equation}
2)Let $N_i$ be a point of a miner increment of sequence $\{{\tilde{a}}(n)\}$ and $p_i$ be the least prime divisor of the product $(N_i-5)(N_i-3).$ If $N_i$ does not complete the run of points of the miner increments after $n_m,$ then the following point of miner increment is
\begin{equation}\label{7.3}
N_{i+1}=\begin{cases}N_i+p_i-1,\enskip if \;\; p_i=2\enskip or\enskip p_i|(N_i-3),
\\N_i+2p_i-5,\enskip  if\;\; p_i>2\enskip and\enskip p_i|(N_i-5).\end{cases}
\end{equation}
3)If the point $N_h$ completes the run of points of miner increments after $n_m,$ then the following point of main increment is \newpage
\begin{equation}\label{7.4}
n_{m+1}=2N_h-4.
\end{equation}
 \end{theorem}
 Note that (7.4) corresponds to (6.5).
 \begin{corollary}
 Conditions of Lemma 3 never satisfy.
 \end{corollary}
 \bfseries Proof. \mdseries  From (7.3) we conclude that after every odd point of miner increment follows even point of miner increment.$\blacksquare$
\begin{remark}
 In connection with Theorem 5 it is interesting to consider a close processes of receiving of twin primes. Let $a$ be odd integer (positive or negative) and $N_i$  be even. Let $p_i$ be the least prime divisor of the product $(N_i-a-2)(N_i-a)$ ( in case of positive $a,\enskip N_i-a-2\geq3). $  Put
$$N_{i+1}=N_i+p-1.$$ One can conjecture that for some $j\geq i,$ the numbers $N_j-a-2,\enskip N_j-a$ will be twin primes. An important shortcoming of such process from the \upshape calculating point of view \slshape is the impossibility to use the formal algorithms for computation of the $\gcd.$
 \end{remark}

 \;\;\;\;\;\;\;\;\;
  \section{Properties of miner increments in supposition of finiteness of twin primes}
 \begin{condition}\label{1}
There exists the maximal second of twin primes $N_{tw}$ such that all seconds of twin primes belong to interval $[5,\enskip N_{tw}].$
 \end{condition}

\begin{corollary}\label{5}
There exists the last point $n_T$ of a main increment of the sequence $\{{\tilde{a}}(n)\}.$
 \end{corollary}
\begin{lemma}\label{4}
If Condition 1 satisfies, then the set of the points righter $n_T$ of nontrivial (miner) increments is infinite.
 \end{lemma}
\bfseries Proof.  \mdseries Suppose that there exists the last point $n=\nu$ of a nontrivial increment, i.e. the set of points of miner the increments is not more than finite. Since  we have
$$\tilde{a}(\nu)=2\nu-4,$$
then for every positive integer $x,$ we find
$$\tilde{a}(\nu+x)=2\nu-4+x.$$\newpage
In particular, for $x=\nu-5,$
$$\tilde{a}(2\nu-5)=3\nu-9.$$
But now the following point $2\nu-4$ is a point of nontrivial increment. Indeed, $\gcd(2\nu-6,\enskip 3\nu-9)=\nu-3.$
Since, evidently, $2\nu-4>\nu,$ then we have contradiction. $\blacksquare$\newline
Besides, from the proof of Lemma 4 the following statement follows.
\begin{lemma}\label{5}
After every $n\geq n_T$ there is not a run of more than $n-5$ trivial increments.
 \end{lemma}
\begin{lemma}\label{6}
Before every nontrivial increment of the magnitude $t$ we have exactly $t-2$ trivial increments.
 \end{lemma}
  \bfseries Proof.  \mdseries Indeed, by the run of formulas (6.2), on every segment
   $$[n_i+l_j+1,\enskip n_i+l_{j+1}-1]$$
we have exactly $l_{j+1}-l_j-1$ points of trivial increments and after that we obtain a nontrivial increment of the magnitude $l_{j+1}-l_j+1.$ $\blacksquare$
\section{Several arithmetical properties of points of the miner increments of sequence $\{{\tilde{a}}(n)\}$}
Further we continue study sequence $\{{\tilde{a}}(n)\}.$
\begin{lemma}\label{7}
 If $M_i$ is an even point of miner increment, then $M_i$ is not multiple of 3.
 \end{lemma}
\bfseries Proof. \mdseries We use induction.
Since $n_m\equiv1\pmod3,$ then, by (8.2), $p_0>3$ and it is easy to see that $M_{1}$ is not multiple of 3. Indeed, in (8.2) it is sufficient to consider cases $p_0\equiv1\pmod3$ and $p_0\equiv2\pmod3.$ Further, using (8.1), note that if the case  $M_i\equiv1\pmod3$ is valid, then the passage from $M_i$ to $M_{i+1}$ is considered as the passage from $n_m$ to $M_1.$ If, finally, $M_i\equiv2\pmod3,$ then $p_i=3,$ and again $M_{i+1}$ is not multiple of 3.$ \blacksquare$
\begin{lemma}\label{8}
 If $N_i$ is an odd point of miner increment, then the congruence $N_i\equiv5\pmod6$ is impossible.
 \end{lemma}
\bfseries Proof. \mdseries Since, by (7.3), after every odd point of miner increment $t$ immediately follows the even point $t+1$ of miner increment, then we should have $N_i+1\equiv0\pmod6.$ This contradicts to Lemma 7.$ \blacksquare$
\begin{lemma}\label{9}
 If $N_i\equiv4\mod6$ is a point of miner increment, then the magnitude of increment in point $N_{i+1}$ is not less than 5.
 \end{lemma}\newpage
 \bfseries Proof. \mdseries Since from Lemmas 7-8 we have $N_{i+1}-N_i\geq3,$ then the lemma follows from Lemma 6.$ \blacksquare$

\begin{lemma}\label{10}
 After every even point of miner increment $N_i$ of the form $N_i\equiv2\pmod6$ follows the odd point $N_i+1$ of miner increment (of the form 6l+3).
 \end{lemma}
 \bfseries Proof. \mdseries Since $N_i-5\equiv0\pmod3,$ then by (7.3), in this case $p_i=3$ and point $N_{i+1}=N_i+2p_i-5=N_i+1$ is the following increment.$ \blacksquare$
 \begin{lemma}\label{11}
The magnitude $\Delta$ of every miner increment either $\Delta=2$ or $\Delta\geq5.$ Moreover, in the second case the previous miner increment has the form $6m+4.$
\end{lemma}
\bfseries Proof. \mdseries From Lemmas 7,8 all points of miner increments have one of the form $6t+i,\enskip i=1,2,3,4.$ Besides, from (7.3) and Lemma 10 the miner increments $\Delta=2$ occur after every points of miner increments of the form $6t+i,\enskip i=1,2,3,$ while, by Lemma 9, after every point of miner increments of the form $6t+4$ we have a miner increment not less than 5.$ \blacksquare$
 \begin{lemma}\label{12}
If Condition 1 satisfies then there are infinitely many points of miner increment of the form $6m+4.$
\end{lemma}
\bfseries Proof. \mdseries  In view of Lemmas 4 and 11, it is sufficient to prove that the process (7.3) which contains  only $p=2$ is finite. Let $N_i$ be point of miner increment 2 such that all follow miner increments are 2. By Lemma 6, it is possible only if all points $N_i, N_i+1,  N_i+2,...$ are points of miner increments. Consider any even point $N_j\equiv1\pmod3,\enskip j\geq i.$ Since $N_j-3$ and $N_j-5$ are not multiple by 2 or 3, then, by (7.3),  $N_{j+1}-N_j>1.$  This contradiction completes the proof.
$ \blacksquare$
\section{Twin numbers of the second kind and accompanying numbers}
\emph{Notation and terminology}. Everywhere below $lpd(n)$ denotes the least prime divisor of $n;$ $p_n$ denotes the $n$-th prime number; $c_i,\enskip i\geq0,$ are constants; $N_{tw}\leq\infty$ is the greater number of the last twin primes pair; $A_1$ is the set of those even $N$ for which $lpd(N-1)<lpd(N-3)$ (cf.A245024 \cite{6}); $A_2$ is the set of those even $N$ for which $lpd(N-1)>lpd(N-3)$ and such that $lpd(N-3), \enskip lpd(N-1)$ are not twin primes (cf.close A243937 \cite{6}); the numbers from the set $A_i$ we call $N_i$-numbers,
 $i=1,2;$ we denote by $N_1(n)$ a $N_1$-number with $lpd(N_1-1)\geq p_n$ and by $N_2(n)$
  a $N_2$-number with $lpd(N_2-3)\geq p_n.$ Finally, we denote by $N_i^{(1)}(n)$ the minimal $N_i(n)$-number,
  \newpage
   $i=1,2$ (cf.A242719, A242720 \cite{6}).\newline
 \indent One can obtain the sequence of twin primes in the following way. Consider
  sequence $\{t_n\}:$"Smallest even $k$ such that $lpd(k-1) > lpd(k-3) \geq p_n, \enskip n\geq2."$ The sequence begins (cf. A242758)
  \begin{equation}\label{10.1}
  6, 8, 14, 14, 20, 20, 32, 32, 32, 44, 44, 44, 62, 62, 62, 62,... \enskip.
  \end{equation}
  Each its term $t(n)$ is associated with a pair of twin primes $t(n)-3, t(n)-1.$ Since the
  lesser numbers of twin primes grow faster than primes, then usually the terms have a multiplicity more than 1. A natural accompanying sequence is "Smallest even $k$ such that $lpd(k-3) > lpd(k-1) \geq p_n, \enskip n\geq2".$ It is sequence $\{N_1^{(1)}(n)\}$ (cf. A242719):
  \begin{equation}\label{10.2}
  10, 26, 50, 170, 170, 362, 362, 842, 842, 1370, 1370,...\enskip.
  \end{equation}
  What sequence could naturally replace sequence $\{t_n\}$ in case $N_{tw}<\infty?$
  Evidently, the sequence "Smallest even $k$ such that the pair $\{k-3,k-1\}$ is not a twin primes pair and $lpd(k-1)>lpd(k-3)>=p_n."$ It is sequence $\{N_2^{(1)}(n)\}$ (cf. A242720):
  \begin{equation}\label{10.3}
  12, 38, 80, 212, 224, 440, 440, 854, 1250, 1460, 1742,...\enskip.
  \end{equation}
  Again a natural accompanying sequence is $\{N_1^{(1)}(n)\}$ (\ref{10.2}).
  The pairs $\{N_2^{(1)}(n)-3, N_2^{(1)}(n)-1\}$ we call \emph{twin numbers of the second kind}. According to (\ref{10.3}), we have the first such pairs
  $$\{9,11\}, \{35,37\}, \{77,79\}, \{209,211\}, \{221,223\}, \{437,439\},$$
  \begin{equation}\label{10.4}
  \{851,853\}, \{1247,1249\}, \{1457,1459\}, \{1739,1741\},...\enskip.
  \end{equation}
\begin{theorem}\label{t6}
There are infinitely many twin numbers of the second kind $\{N_2^{(1)}(n)\}$ and infinitely
many their accompanying numbers $\{N_2^{(1)}(n)\}.$
 \end{theorem}
\begin{proof}
 It is sufficient
  to prove the infinity of $N_i(n),\enskip i=1,2.$
  Consider sequences
  \begin{equation}\label{10.5}
  x_n=(p_n-1)!+2, \enskip y_n=(p_n-2)!+2, \enskip n\geq3.
  \end{equation}
  By the Wilson theorem and one of its corollary, we have
  $$x_n-1=(p_n-1)!+1\equiv1 \pmod {p_i},\enskip i\leq n-1,\enskip and\enskip \equiv0 \pmod{p_n},$$
  $$x_n-3=(p_n-1)!-1\equiv-1 \pmod {p_i},\enskip i\leq n.$$
  So, $$lpd(x_n-3)>lpd(x_n-1)=p_n$$
  and we conclude that $x_n$ is a $N_1$-number.
  Analogously,
   $$y_n-1=(p_n-2)!+1\equiv1 \pmod {p_i},\enskip i\leq n,$$
    \newpage
  $$y_n-3=(p_n-2)!-1\equiv-1 \pmod {p_i},\enskip i\leq n-1, \enskip and\enskip \equiv0 \pmod{p_n}.$$
  So,
  $$lpd(y_n-1)>lpd(y_n-3)=p_n$$
  and  $y_n$ is a $N_2$-number.
  \end{proof}
  The second proof is based on Dirichlet theorem on arithmetical progressions (cf. comment
  by R. Israel in A242033\cite{6}).
  \begin{proof}
  Consider the congruence $p_nx\equiv1 \pmod {\prod_{i<n-1}p_i}, \enskip n\geq2.$
Let $q>p_n$ be a prime solution which exists by Dirichlet's theorem. Now $p_nq-2$ is
divisible by none of primes less than or equal $p_n.$
Hence, $p_nq+1$ is a $N_1(n)$-number. Indeed, $p_n=lpd((p_nq+1)-1)<lpd(p_nq-2).$
\newline
 \indent Further, consider the congruence $p_nx\equiv-1 \pmod {\prod_{i<n-1}p_i}$ and let $r>p_n$ be a prime solution. Now $p_nr+2$ is divisible by none of primes less than or equal $p_n.$
Hence, $p_nr+3$ is a $N_2(n)$-number. Indeed, $p_n=lpd((p_nr+3)-3)<lpd(p_nr+2).$
\end{proof}

  In addition note that, by the definition, we have
   \begin{equation}\label{10.6}
  N_i^{(1)}(n)\geq p_n^2+1,\enskip i=1,2.
  \end{equation}
  The equality satisfies for $N_1^{(1)}(n),$ if $p_n^2-2$ is prime.\newline
\indent In the second part of the paper, the important role plays the following our postulate.
\begin{postulate}\label{Pos1}
At least for infinite set of $n,$ we have
\begin{equation}\label{10.7}\
\max(N_1^{(1)}(n),\enskip N_2^{(1)}(n))<(\min(N_1^{(1)}(n),\enskip N_2^{(1)}(n)))^2.
\end{equation}
\end{postulate}
Let us indicate a sufficient condition for truth of the Postulate.

\begin{theorem}\label{t7}
 If, for arbitrary odd prime $P=p_n,$ there exist primes $Q, R$ in interval $[P, P^3]$
  such that both numbers $PQ+2$ and $PR-2$ are primes, then the Postulate holds.
  \end{theorem}
  \begin{proof}
   For $P=p_n,$  $N_1(n)$-numbers, by the definition,
 possess property $lpd(N_1(n)-3)\geq lpd(N_1(n)-1\geq P,$ while for $N_2(n)$-numbers we have
 $lpd(N_1(n)-1)\geq lpd(N_1(n)-3\geq P,$ and, by (\ref{10.6}), every $N_1, N_2\geq P^2+1.$ Since,
 by condition, $PQ+2$  is prime, then $PQ+3$ is $N_2(n)$-number; since $PR-2$ is prime, then $PR+1$ is $N_1(n)$-number. Besides, for our $N(n)$-numbers, we have
$$P^2+1\leq PQ+3=N_2(n)\leq P^4+3<(P^2+1)^2$$
and
 \newpage
 $$P^2+1\leq PR+1=N_1(n)\leq P^4+1< (P^2+1)^2.$$
Thus also both the minimal $N_1(n)=N_1^{(1)}(n)$ and the minimal $N_2(n)=N_2^{(1)}(n)$ are in interval $[P^2+1, (P^2+1)^2].$
Let, say, $N_1^{(1)}(n)<N_2^{(1)}(n).$ Then
$$N_2^{(1)}(n)\leq (P^2+1)^2\leq (N_1^{(1)}(n))^2.$$
But easily to show that $N_2^{(1)}(n)\neq (N_1^{(1)}(n))^2,$ i.e., $PR+1\neq(PQ+3)^2,$ or
$R=PQ^2+6Q+8/P>P^3,$ which contradicts the condition. So $N_2^{(1)}(n)<(N_1^{(1)}(n))^2.$
\end{proof}

 \section{A heuristic proof of the postulate for large $n$}
  Consider a progression
\begin{equation}\label{11.1}
 F(p_n,t)=2p_nt-(p_n-2)=p_n(2t-1)+2, \enskip t=1,2,...\enskip .
\end{equation}
The number of primes of such a form not exceeding $y$ is
\begin{equation}\label{11.2}
\sim (p_n-1)^{-1}y/\ln y \enskip\enskip (y\rightarrow\infty).
\end{equation}
Formally, the probability for $F$ to be prime grows with the number of prime divisor
of $2t-1.$ Therefore, $F$ is prime more often when $2t-1$ is composite number, than
it is prime. It is well known that the number $\omega(2t-1)$ of prime divisors of $2t-1$
in average is $\ln\ln (2t-1).$
 Since $F\leq y$ yields $2t-1\leq\frac{y-2}{p_n}$ and $2t-1$ runs all odd integers in the interval $(0,\frac{y-2}{p_n}],$
then $2t-1$ runs all primes in this interval. So, "primarity coefficient" of $2t-1,$ when
$F\leq y,$ is \enskip $2\pi ((\frac{y-2}{p_n})/(\frac{y-2}{p_n}))$ and, if do not take into account the noted dependence of primarity of $F$ from the number of prime divisors of
 $2t-1,$ then the number $E(y)$ of primes $F(p_n,t)\leq y$ with primes $2t-1$ would be
$$E(y)\sim 2\pi ((\frac{y-2}{p_n})/(\frac{y-2}{p_n}) (p_n-1)^{-1}y/\ln y$$
 \begin{equation}\label{11.3}
\sim \frac{2y}{(p_n-1)(\ln y)^2}.
 \end{equation}
  But, taking into account this factor, we can suppose that it acts proportionally
 to $\omega(2t-1).$ Besides, the record values of $\omega(2t-1)$ arise when $2t-1$
  is the product of the first several consecutive odd primes. In this case we have \cite{9}
   $\omega(2t-1)\sim \ln(2t-1)/\ln\ln(2t-1).$
  So, instead of (\ref{11.3}), it is natural to
 expect that at least the following inequality hods
\begin{equation}\label{11.4}
E(y)\geq \frac{c_{0}y\ln\ln y}{(P-1)(\ln y)^3}.
\end{equation}
Now we set $y=P^4.$ Since now $2t-1$ runs all odd integers in the interval $(0,2\frac{P^4-2}{P}],$ then we can choose from this interval a prime $q\geq p_n$ such that
\newpage
\begin{equation}\label{11.5}
 F=p_nq+2< p_n^4
  \end{equation}
 is prime. This means that $F+1$ is a $N_2(n)$-number.
 \begin{remark}\label{r9}
 Linnik \cite{11}-\cite{12} proved that the least prime $p(a,d)$ in the progression $a+dt$
 does not exceed $Cd^L,$ where $C,L$ are absolute constants. Without GRH
 Triantafyllos \cite{13} proved only that $L=5.$ It is the best result without
 GRH. Using this result, we cannot guarantee the existence of prime $F=p_nq+2$ which is now less than $p_n^4.$ But, using GRH, Heath-Brown \cite{10} proved that
  \begin{equation}\label{11.6}
  p(a,d)\leq (1+o(1))(\varphi(d)\ln d)^2,
 \end{equation}
 where $\varphi$ is the Euler totient function.
 In our case this means that $F$ could be chosen in interval
 $c_1((p_n-1)\ln p_n)^2\leq F< p_n^4.$
 \end{remark}
\indent Furthermore, by the analogous arguments, considering  a progression
 \begin{equation}\label{11.7}
 F_1(p_n,t)=2p_nt-(p_n+2), t=1,2,...\enskip ,
\end{equation}
we find a prime $r\geq p_n$ such that
\begin{equation}\label{11.8}
 F_1=p_nr-2< p_n^4,
\end{equation}
 is prime and, consequently, $F_1+3$ is a $N_1(n)$-number.
 Thus also both the
 minimal $N_1(n)=N_1^{(1)}(n)$ and the minimal $N_2(n)=N_2^{(1)}(n)$ are in interval
 $[p_n^2+1, (p_n^2+1)^2]$ and we have either $N_1^{(1)}(n)<N_2^{(1)}(n)<(N_1^{(1)}(n))^2$ or $N_2^{(1)}(n)<N_1^{(1)}(n)<(N_2^{(1)}(n))^2.$ $\blacksquare$\newline

 \section{ Tolev's theorem}
 In 1999, Tolev \cite{7} proved the following theorem.
  \begin{theorem}\label{t8} (\cite{7})
 For a constant $c_0>0,$ there are at least $c_0x^2/(\ln x)^6$ triples of primes
 $\{q_1,q_2,q_3\}$ in interval $(x,2x),$ satisfying $q_1+q_2=2q_3$ and such that
 $\min (lpd (q_1+2), lpd (q_2+2))\geq x^{0,167}$ and $lpd (q_3+2)\geq x^{0.116}.$
 \end{theorem}
  Note that Theorem \ref{t8} is based on a lower estimate $(x^2/(\ln x)^3)$ of a generalized Chebyshev's function
$$\Gamma=\sum \ln p_1 \ln p_2 \ln p_3,$$
where the summing is over $x < p_1, p_2, p_3 < 2x$ such that $p_1 + p_2 = 2p_3$ and, if $z_i = x^{\alpha_i}$ , where $\alpha_i, \enskip i = 1, 2, 3$ are some constants from the interval $(0, 1/4),$ then
$p_i + 2$ is divisible by none of odd primes less than $z_i, \enskip i=1, 2, 3$.
Reading the proof of Theorem \ref{t8} \cite{7}, one can see that it does not depend on the changing $p_i + 2$ by $p_i - 2.$ So, the following symmetrical theorems hold.
\newpage
\begin{theorem}\label{t9}
For a constant $c_3>0,$ there are at least $c_0x^2/(\ln x)^6$ triples of primes
 $\{q_1,q_2,q_3\}$ in interval $(x,2x),$ satisfying $q_1+q_2=2q_3$ and such that
 $\min (lpd (q_1-2), lpd (q_2-2))\geq x^{0,167}$ and $lpd (q_3+2)\geq x^{0.116}.$
 \end{theorem}
 \begin{theorem}\label{t10}
For a constant $c_4>0,$ there are at least $c_6x^2/(\ln x)^6$ triples of primes
 $\{q_1,q_2,q_3\}$ in interval $(x,2x),$ satisfying $q_1+q_2=2q_3$ and such that
 $\min (lpd (q_1+2), lpd (q_2-2))\geq x^{0,167}$ and $lpd (q_3+2)\geq x^{0.116}.$
 \end{theorem}

\section{An estimate for $N_i(n)$-numbers, i=1,2, in case $N_{tw}<\infty$ }
  Note that, every $q_1$ and $q_2$ in Theorem \ref{t8}, evidently, cannot run less than
 $c_2^{1/2}x/(\ln x)^3$ different values. So, the number of different
 values of $q_1$ in interval $(x,2x)$ is $\geq c_2^{1/2}x/(\ln x)^3.$ \newline
 \indent Set now $x=x(n)=p^{5.989}_{n}.$ Then $p^{5.989}_{n} <q_1<2p^{5.989}_{n}.$
  According to Theorem \ref{t8}, we have
 $$lpd (q_1+2)\geq x^{0.167}=(p^{5.989}_{n})^{0.167}=p_n^{1.000163}>p_n.$$ This yields that every $q_1+3$  is $N_1(n)$-number. Indeed, $lpd((q_1+3)-3)=q_1>
 lpd((q_1+3)-1)>p_n.$ Thus $N_1(n)\leq q_1+3<2x+3<p_n^6.$
 Analogously, using Theorem \ref{t9} for $x=x(n)=p^{5.989}_{n}$ and noting that in this theorem
in case, when $q_1-2, q_1$
 are not twin primes, every $q_1 + 1$ is $N_2(n)$-number, we obtain $N_2(n)<p_n^6.$\newline
 \indent Thus, if $N_{tw}<\infty,$ and $p_n>N_{tw},$ then we have
 \begin{equation}\label{13.1}
N_i^{(1)}(n)<p_n^6,\enskip i=1,2.
\end{equation}

\section{A statistical symmetry between $N_1(n)$ and $N_2(n)$-numbers}
Let $N$ be positive even number such that
\begin{equation}\label{14.1}
N\equiv a_2b_2\frac{M}{3}+...+a_{n-1}b_{n-1}\frac{M}{p_{n-1}}+a_nb_{n}\frac{M}
{p_{n}} \pmod {M},
\end{equation}
where $ M=M_n=\prod_{i=1}^{n} p_i,\enskip , \enskip b_i\frac{M}{p_i}\equiv1
\pmod{p_i}$ and integers $a_i$ are nonnegative residue modulo $p_i$ respectively,
 such that $a_i\neq1,3 \pmod{p_i},\enskip i=2,...,n-1,$ while $a_n$ is an arbitrary
nonnegative residue modulo $p_n.$ By Chinese theorem, the least prime divisors of both numbers $N-1$ and $N-3$ $(lpd(N-1)$ and $lpd(N-3))$ are equal or more than $p_n.$
Thus, according to our notation, $N$ is
 $N(n)$-number.
 Consider firstly the case $N_{tw}<\infty.$ Let $n$ be such that
 \begin{equation}\label{14.2}
 p_{n}>N_{tw}.
  \end{equation}
Evidently,
\begin{equation}\label{14.3}
 N(n)\in (p^2_n, M_n].
  \end{equation}
  \newpage
 The number $m_n$ of all different considered $N(n)$-numbers is
 \begin{equation}\label{14.4}
  m_n=(p_2-2)(p_3-2)...(p_{n-1}-2)p_n.
  \end{equation}
  Moreover, by the symmetry with respect to $N(n)-2,$
   we have approximately the same number of $N(n)$-numbers for which
   $lpd(N(n)-3)>lpd(N(n)-1)$ and of $N(n)$-numbers for which  $lpd(N(n)-3)<lpd(N(n)-1),$ and
 these types of $N(n)$-numbers, i.e., $N_1(n)$ and $N_2(n)$ have approximately the same distribution.
   \begin{remark}\label{r5}
  This symmetry manifests itself stronger especially in the situation when,
  by the condition $(\ref{14.2}),$ in the interval $(p^2_n, M_n] $
 there are no twin primes. Indeed, if $(N-3, N-1)$ is a pair of twin primes, then
 a priori we have $lpd(N-1)>lpd(N-3).$ However, if to write $N-3 \enskip'=\enskip'\enskip N-1$
 (and only for them) and to include also $'=\enskip' $ in the definition of $N$-numbers,
 i.e., to include $N$-numbers with $'=\enskip'$ in both types of $N$-numbers,
 then even for small $n$, for example, in case $n=4,\enskip p_4=7,$ considering the interval
 $(49, 210],$ we obtain the following $N$-numbers:
    $\{50,62,74,80,92,104,110,122,134,140,152,164,170,182,194,200\}.$ It is interesting that the $N$-numbers with strong inequalities $lpd(N-1)<lpd(N-3)$ and  $lpd(N-1)>lpd(N-3)$ here alternate. See also sequences $A243803, A243804$ and especially $A242974$ \cite{6}.
    \end{remark}
 Since the average distance $\rho(n)$ between two consecutive $N(n)$-numbers
in interval $( p^2_n, M_n]$ is not more than $\frac{M_n}{m_n},$ then we have
\begin{equation}\label{14.5}
\rho(n)\leq \frac{ M_n}{m_n}=2\prod_{i=2}^{n-1}(1+\frac{2}{p_i-2})\leq5.2826...\prod_{i=2}^{n-1}(1+\frac{2}{p_i}),
  \end{equation}
since
$$ \prod_{i=2}^{n-1}(1+\frac{2}{p_i-2})/\prod_{i=2}^{n-1}(1+\frac{2}{p_i})=
  \prod_{i=2}^{n-1}(1+\frac {4}{p_i^2-4})$$
\begin{equation}\label{14.6}
  <\prod_{i=2}^{\infty}(1+\frac {4}{p_i^2-4})=2.6413...\enskip .
\end{equation}
Furthermore, by a Rosser result \cite{1}, we have
\begin{equation}\label{14.7}
\prod_{i=2}^{n-1}(1-\frac{2}{p_i})=\frac{0.832429...+o(1)}{\ln^2p_{n-1}}.
\end{equation}
Besides,
\begin{equation}\label{14.8}
\prod_{i=2}^{n-1}(1-\frac{2}{p_i})\prod_{i=2}^{n-1}(1+\frac{2}{p_i})=C+o(1),
\end{equation}
where
\newpage $$C=\prod_{i=2}^{\infty}(1-\frac{4}{p^2_i})=0.3785994... $$

and, by (\ref{14.7})-(\ref{14.8}) (we have here a very large $n$) we find
$$\prod_{i=2}^{n-1}(1+\frac{2}{p_i})=\frac{0.3785994...+o(1)}{0.832429...+o(1)}
\ln^2 p_{n-1}$$
$$\leq 0.4549\ln^2 p_{n-1}. $$
Thus, according to (\ref{14.5}), we have
\begin{equation}\label{14.9}
 \rho(n)\leq\frac{ M_n}{m_n}\leq2.4026\ln^2 p_{n-1}.
\end{equation}
\indent In case of the infinity of twin primes, the average distance between them on the interval
$(p_n^2, M_n)$ is more than $C\ln^2 M_n\gg\ln^2 p_{n-1}.$ For large $n,$ it counts
by already made rounding the result.

\section{A theorem}
\begin{theorem}\label{t11}
Let $N_{tw}<\infty.$ If the Postulate does not satisfy, then it is possible only a finite number of changing
of sign of the difference $d_n=N_1^{(1)}(n)-N_2^{(1)}(n).$
\end{theorem}
\begin{proof}
By the condition, there exist $n_1$ such that, for $n\geq n_1,$ the Postulate does
not satisfy. Suppose that after $n_1$ we have a change of sign of $d_n.$  Consider two consecutive numbers $n-1$ and $n$ such that $p_{n-1}^2<N^{(1)}_2(n-1)<N^{(1)}_1(n-1)$ and
 $p_n<N^{(1)}_1(n)<N^{(1)}_2(n).$ Since the Postulate does not satisfy, from the first inequality we have
 $$p_{n-1}^4<(N_2^{(1)}(n-1))^2<N_1^{(1)}(n-1).$$
 From the second inequality we have
 $$N_2^{(1)}(n)>(N_1^{(1)}(n))^2.$$
 In view of $N_1^{(1)}(n)$ and $N_2^{(1)}(n)$ are nondecreasing, then further we
 have
 $$N_2^{(1)}(n)>(N_1^{(1)}(n))^2\geq $$
 $$(N_1^{(1)}(n-1))^2>(N_2^{(1)}(n-1))^4\geq p^8_{n-1}$$
 However, for $p_n>N_{tw}$ this contradicts the estimate (\ref{13.1}).
\end{proof}
\newpage
\section{AiB-axiom}
Suppose that we have two unprovable but very plausible conjectures $A$ and $B$.
There is a sense to accept also an unprovable conjecture that $A$ implies $B$ as an axiom
(we call it an axiom of type "AiB"), if it leads to a consistent meaningful theory, such that
in its frameworks we prove that also $(\overline{A}\Rightarrow {B}).$ Thus, by the AiB-axiom,
$A$ is a sufficient condition for $B$ and, if this sufficient condition does not satisfy,
then $B$ also takes place. \newline
\indent In our case, $A$ is a very plausible inequality (\ref{10.7}) which we call
"postulate", and $B$ is "the infiniteness of twin primes".
\begin{remark}\label{r5}
I found an error in proof of the former Theorem $6 (2010)$ stating that $"A\Rightarrow B".$
Since, despite my best potential efforts to correct it, I was not able
to find a right proof, I began to consider this error as unrecoverable one. However,
 "Theorem 6" led me, using Chinese and remarkable Tolev's theorems, to an interesting
 theory, including reducing the supposition of the finiteness of twin primes to an arbitrary long coin-flipping experiment in which only "heads" appear (see version $34$ of this paper, where Theorem 6 should be replaced by the considered axiom; note that I mean, namely this statement, when I say that also $\overline{A}$ yields $B$). So I naturally became to idea of "AiB-axiom".
\end{remark}
\begin{remark}\label{r6}
Consider an example of connection between $N_i^{(1)}(n),\enskip 1,2,$ and twin primes.
Note that $p_n^2+1$ is $N_1^{(1)}(n)$-number, if and only if $p_n^2-2$ is prime.
Let, furthermore, for $k\leq p_n+1,\enskip k^2+2$ be $N_2^{(1)}(n)$-number, $n\geq3.$
Then $k-1, k+1$ are twin primes and also $k^2+1$ is prime. Indeed, by the 
definition, $k^2+2=N_2^{(1)}(n),$ if and only if $k$ is the minimal with the 
condition $lpd((k^2+2)-1)>lpd((k^2+2)-3)\geq p_n.$ Then such $k$ is unique, such that
$k-1=p_n,\enskip k+1=p_{n+1}.$ Moreover, since $lpd(k^2+1)\geq k+1,$ then $k^2+1$ is prime. 
 Such suitable values of $k$
are $($cf. $A070155$ \cite{6}$)$
\begin{equation}\label{16.1}
6,150,180,240,270,420,570,1290,1320,...\enskip.
\end{equation}
\indent In connection with sequence $A070155,$ note that the case $n=2,\enskip
p_n=3,\enskip k=4=A070155(1),$ when $N_2^{(1)}(n)=p_n^2+3=12<18,$ is a special, since, for $n>2,$ $p_n^2+2\equiv0 \pmod3.$\newline
In order to have the considered $N_1^{(1)}$ and $N_2^{(1)}$ in the same values of
$n,$ we should require $p_n^2-2=(k-1)^2-2$ to be prime. Then we obtain the following
 sequence, instead of $(\ref{16.1}):$
\begin{equation}\label{16.2}
6,240,570,1290,2310,2550,2730,3360,...\enskip.
\end{equation}
\newpage
Thus, if this sequence is infinite, then $(\ref{10.7})$ satisfies together with the infiniteness of twin primes. Construction of this sequence is a some additional "motivation" of the axiom of type AiB.\newline
\end{remark}
\section{Conclusive remarks and problems}
 It is highly interesting that for numbers $a(n)=N^{(1)}_1(n)=A242719(n),$
$b(n)=N^{(1)}_2(n)=A242720(n),$ most likely, it follows that
\begin{conjecture}\label{con4}
For $n\geq2, \enskip a(n)-3$ is prime and $a(n)-1$ is semiprime;
 for $n\geq21, \enskip b(n)-3$ is semiprime and $b(n)-1$ is prime.
\end{conjecture}
Thus, especially, sequence A242719 is a beautiful illustration of the very
known Chen's result \cite{8} in this direction. Chen proved that there exist infinitely many primes $p$ such that $p+2$ is prime or semiprime.\newline
\indent Note that, Conjecture \ref{con4} was verified by J. C. Moses
up to 2001 and, respectively, up to 2501 for $a(n)$ and, respectively, for $b(n).$ Before $n=2501,$
he found only two semiprimes of the form $b(n)-1:$ $b(16)-1 = 4189 = 59\cdot71$ and $b(20)-1 = 6889 = 83^2.$\newline
\indent In connection with Conjecture \ref{con4}, let us show how to find $lpd(a(n)-1)$
and $lpd(b(n)-3).$ With this aim, consider sequence $\{\alpha(n)\}, \enskip n\geq2,$
such that $\alpha(n)$ is the smallest even $k$ for which $lpd(k-1)=p_n,$ while $lpd(k-3) >p_n$ (cf.A242489\cite{6}). Passing from this non-monotonic sequence to the nondecreasing sequence $A242719=\{a(n)\},$ we notice that $\{a(n)\}$ consists of chains of different lengths $s\geq1,$ such that each chain consists of the same numbers
$a(k)=a(k+1)=...=a(k+s-1).$ The last term of the chain $a(k+s-1)=\alpha(k+s-1)$ is a term of $\{\alpha(n)\}=A242489$ and, therefore, is divisible by $p_{k+s-1}.$ Note that $lpd(a(k+s-1)-1=p_{k+s-1}.$
Thus, in order to find $lpd(a(n)-1$ over A242719 we should find the last term $a(m)=a(n)$ of the chain which contains $a(n).$ Now $lpd(a(n)-1)=p_m.$ Analogously we find $lpd(b(n)-3)$ over A242720 (cf. A242490). By the way, we conjecture that in each sequences A242719, A242720 there are arbitrary long such chains.
\indent Finally, instead of (\ref{13.1}), we conjecture that
\begin{equation}\label{17.1}
\max(a(n),\enskip b(n))< p_n^4, \enskip n\geq2.
\end{equation}
Moreover, there are bases to think (cf. Remark 6) that
\begin{equation}\label{17.2}
 \max(a(n),\enskip b(n))=O(n^2(\log n)^2).
\end{equation}
\newpage
 \section{To the reader}
I apologize that I did so many versions of the paper. I worked step by step, since
in my current situation I cannot leave "on then" unfinished thoughts. Sometimes, I did
 stupid mistakes and should was correct them, increasing the number of versions.
  However, while working on this paper, I received really a great fun and I hope that it at least a little was transmitted to the reader.

\section{Acknowledgment}
The author sincerely thanks N. J. A. Sloane for his best edition of author's
sequences connected with this search. He is also grateful to Peter J. C.
Moses (England) for his important permanent computer help.

\end{document}